\documentclass[12pt]{article}
\usepackage[margin=2cm]{geometry}
\usepackage{amsfonts}
\usepackage{amssymb}
\usepackage{amsthm}
\usepackage{amsmath}
\usepackage{latexsym}
\usepackage{color}

\begin{document}
\newtheorem{lemma}{Lemma}
\newtheorem{pron}{Proposition}
\newtheorem{thm}{Theorem}
\newtheorem{Corol}{Corollary}
\newtheorem{exam}{Example}
\newtheorem{defin}{Definition}
\newtheorem{re}{Remark}
\newcommand{\la}{\frac{1}{\lambda}}
\newcommand{\sectemul}{\arabic{section}}
\renewcommand{\theequation}{\sectemul.\arabic{equation}}
\renewcommand{\thepron}{\sectemul.\arabic{pron}}
\renewcommand{\thelemma}{\sectemul.\arabic{lemma}}
\renewcommand{\thethm}{\sectemul.\arabic{thm}}
\renewcommand{\there}{\sectemul.\arabic{re}}
\renewcommand{\theCorol}{\sectemul.\arabic{Corol}}
\renewcommand{\theexam}{\sectemul.\arabic{exam}}
\renewcommand{\thedefin}{\sectemul.\arabic{defin}}
\def\REF#1{\par\hangindent\parindent\indent\llap{#1\enspace}\ignorespaces}
\def\lo{\left}
\def\ro{\right}
\def\be{\begin{equation}}
\def\ee{\end{equation}}
\def\beq{\begin{eqnarray*}}
\def\eeq{\end{eqnarray*}}
\def\bea{\begin{eqnarray}}
\def\eea{\end{eqnarray}}
\def\d{\Delta_T}
\def\r{random walk}
\def\o{\overline}
\def\N{\mathbb{N}}
\def\P{\mathbb{P}}
\def\J{\mathcal{J}}

\title{\large\bf On the structure of a class of distributions obeying the principle of a single big jump
\thanks{Research supported by National Science Foundation of China
(No.11071182 ), Natural Science Foundation of the Jiangsu Higher
Education Institutions of China (No.13KJB110025).}}
\author{\small
Hui Xu $^{1}$,~Michael Scheutzow $^{2}$~Yuebao Wang$^{1}$\thanks{Corresponding author.
Telephone: 86 512 67422726. Fax: 86 512 65112637. E-mail:
ybwang@suda.edu.cn},~Zhaolei Cui $^{3}$
\\
{\footnotesize\it 1. School of Mathematical Sciences, Soochow
University, Suzhou 215006, China}\\
{\footnotesize\it 2. Institut f\"ur Mathematik, Technische Universit\"at Berlin, 10623 Berlin, Germany}\\ {\footnotesize\it 3. Changshu Institute of Technology, Changshu 215500, China}}
\date{}

\maketitle
{\noindent\small {\bf Abstract }}

In this paper, we present several heavy-tailed distributions belonging to the new class $\J$ of distributions obeying the principle of a single big jump introduced by Beck et al.~\cite{BBS13}.
We describe the structure of this class from different angles. First, we show that heavy-tailed distributions in the class $\J$ are automatically {\em strongly heavy-tailed} and thus have tails which are not too irregular. Second,
we show that such distributions are not necessarily weakly tail equivalent to a subexponential distribution. We also show that the class of heavy-tailed distributions in $\J$ which are neither long-tailed nor dominatedly-varying-tailed is not only non-empty but even quite rich in the sense that it has a nonempty intersection with several other well-established classes. In addition, the integrated tail distribution of some particular of these distributions shows that the Pakes-Veraverbeke-Embrechts Theorem for the class $\J$ in \cite{BBS13} does not hold trivially.\\



\section{Introduction}
\setcounter{equation}{0}\setcounter{Corol}{0}\setcounter{lemma}{0}\setcounter{re}{0}
In this paper, all distributions have unbounded support contained in $[0,\infty)$.
Recall that a distribution $F$ is called {\em heavy-tailed}, denoted by $F\in \mathcal{K}$, if for all $\alpha>0$,
\begin{eqnarray*}
\int_{0}^{\infty}e^{\alpha y}dF(y)=\infty;
\end{eqnarray*}
otherwise, $F$ is called {\em light-tailed}, denoted by $F\in \mathcal{K}^c$.
Recently,
Beck et al.~introduced the following new distribution class $\J$ in \cite{BBS13}.

Let $\{X_i,i\geq 1\}$ be a sequence of independent and identically distributed (i.i.d.) random variables with common distribution $F$.
Define the class $\J$ as the set of those distributions $F$ such that for all $n\geq 2$,
\begin{eqnarray}\label{def-e101}
\lim_{K\to\infty}\liminf\P\big(X_{n,1}>x-K \mid \sum_{i=1}^n X_i>x\big)=1,
\end{eqnarray}
where $X_{n,k}$ means the $k$-th largest random variable in the sequence $\{X_i,1\leq i\leq n\},1\leq k\leq n$. Here and in the following all
unspecified limits are to be understood for $x \to \infty$.
Beck et al.~\cite{BBS13} show that each of the following two properties are equivalent to (\ref{def-e101}):
\begin{align}\label{def-e104}
&\lim_{K\to\infty}\liminf\P(X_{n,2}\leq K \mid \sum_{i=1}^n X_i>x)=1,\\
&\lim\P (X_{n,2}>g(x)\mid \sum_{i=1}^n X_i>x)=0
\mbox{ for all } g\nearrow \infty,\nonumber
\end{align}
for any (and hence for all) $n\ge 2$.

Note that this definition is
a way of stating formally that the distribution of $X_1$ obeys the {\em principle of a single big jump} which means that conditional on
the sum $X_1+...+X_n$ being unusually large, the probability that a single summand dominates the sum is close to one while the conditional
law of the second largest summand remains tight (the letter $\J$ stands for {\em jump}).

It is natural to ask how the class $\J$ is related to well-established classes like $\mathcal{S}$, $\mathcal{D}$ and $\mathcal{L}$ which are respectively called
{\em subexponential}, {\em dominatedly varying} and {\em long tailed} (for their definition, see \cite{EKM97}). In addition, the class
$\cal{OS}$ of {\em generalized subexponential} distributions first introduced in \cite{K90} is of interest. By definition, it consists of those distributions for which
\begin{eqnarray}\label{def-e106}
C^*(F):=\limsup\overline{F^{*2}}(x)(\overline{F}(x))^{-1}<\infty,
\end{eqnarray}
where $F^{*n}$ denotes the $n$-fold convolution of $F$ with itself for $n\ge2$, and $\overline{F}:=1-F$ denotes the tail of $F$. Note that the class $\mathcal{S}$ corresponds to the
case $C^*(F)=2$.

The following relations for these classes are known:
\begin{equation}\label{inclusion}
\mathcal{S} \subset \mathcal{L} \subset \mathcal{K},\; \mathcal{S} \cup \mathcal{D} \subset \mathcal{J} \subset \cal{OS},\; \mathcal{S}=\mathcal{J} \cap \mathcal{L},\; \mathcal{D}\subset \mathcal{K}
\end{equation}
(see  \cite{EKM97} respectively \cite{BBS13} for relations not involving $\J$ respectively those involving $\J$). It is not true that all distributions in $\J$ are
heavy-tailed (see \cite{BBS13}). The recent paper  \cite{XSW14} actually shows that the class of light-tailed distributions in $\J$ is considerably larger than the union of
the well-known classes $\mathcal{S}(\gamma)$.

In this paper, the object of study is the class $\J\cap\cal{K}$. The list of relations \eqref{inclusion} says that
$\mathcal{J}\cap \mathcal{K}$ contains $\mathcal{S} \cup \mathcal{D}$. We will not only  show that this inclusion is proper but that it is even quite large thus suggesting
that the class $\J$ cannot be simply characterized via other known classes. On the other hand, we will show first that $\mathcal{J}\cap \mathcal{K}$ is (strictly)
contained in the class of {\em strongly heavy-tailed} distributions which are characterized by the property that for all $\lambda>0$ we have
\begin{eqnarray}\label{108}
\lim e^{\lambda x}\overline{F}(x)=\infty.
\end{eqnarray}
This class was denoted by $\mathcal{K}^*$ in \cite{BBS13} and by  $\mathcal{DK}^c$ in  \cite{WCY11}. It is clearly contained in the class of
heavy-tailed distributions but excludes some members of $\mathcal{K}$ with irregular tails. The following result will be proved in Section 2.

\begin{thm}\label{thm101}
The following inclusions hold:
$$\J\cap\mathcal{K}\subset\mathcal{OS}\cap\mathcal{K} \subset \mathcal{DK}^c.$$
\end{thm}
Note that the first inclusion is clear from \eqref{inclusion}, so we only need to show the second one. In addition, we will show in Example \ref{exam300}
that this inclusion is proper.  We remark that
Theorem \ref{thm101} shows in particular that condition (iii)  $F^I \in \J\cap\mathcal{DK}^c$  in Theorem 19 in \cite{BBS13} can be replaced by the equivalent condition
$F^I \in \J\cap\mathcal{K}$.\\

Before we state our second result, we define the class $\mathcal{DK}_1$ which was introduced by  Wang et al.~\cite{WCY11}, as the set of all distributions which
satisfy
\begin{eqnarray}\label{109}
\lim x^{\delta}\overline{F}(x)=\infty,\mbox{ for some } \delta>0.
\end{eqnarray}
Note that $\mathcal{D} \subset \mathcal{DK}_1 \subset \mathcal{DK}^c$. We will call two distributions $F$ and $G$ {\em weakly tail equivalent}, denoted by
$\overline{F}\thickapprox\overline{G}$, if  $$0<\liminf\overline{F}(x)(\overline{G}(x))^{-1}\le\limsup\overline{F}(x)(\overline{G}(x))^{-1}<\infty.$$

\begin{thm}\label{thm102}
The class $(\J\cap\cal{K})\setminus(\cal{L}\cup\cal{D})$ is non-empty.  Moreover, none of its intersections with  $\mathcal{DK}_1$ and its complement is empty and
each of these two subclasses contains both distributions which are weakly tail equivalent to a distribution in $\mathcal{S}$ and distributions which are not.
\end{thm}

We will provide four corresponding examples in Sections 3 and 4. Note that it does not matter whether or not we replace $\mathcal{L}$ by $\mathcal{S}$
in Theorem \ref{thm102} since $\mathcal{S}=\J \cap \mathcal{L}$ by \eqref{inclusion}.\\

Finally, we investigate the class $\J$  with respect to integrated tail distributions. In
Theorem 19 of \cite{BBS13}, the integrated tail distribution of the claim size in the Sparre Andersen risk model is required to belong to the class $\J\cap\cal{K}$,
so the question arises whether there exists a distribution $F$ whose integrated tail distribution $F^I\in \J\cap\cal{K}\setminus\cal{S}$?  Otherwise, if $F^I\in \cal{S}$,
then the corresponding result is the known Pakes-Veraverbeke-Embrechts Theorem, see Theorem 16 of \cite{BBS13}.
\vspace{0.2cm}

To answer this question, we recall the concepts of an integrated tail distribution and a generalized long-tailed distribution.

For a distribution $F$, if $0<\mu:=\int_0^\infty \overline{F}(y)dy<\infty$, then the distribution $F^{I}$ defined by
$$\overline{F^{I}}(x)={\mu}^{-1}\int_x^{\infty}\overline{F}(y)dy\emph{\textbf{\emph{1}}}(x\ge 0), ~x\in [0,\infty)$$
is called the {\em integrated tail distribution} of $F$.

A distribution $F$ is called {\em generalized long-tailed}, denoted
by $F\in{\cal{OL}}$ (see \cite{SW05}), if for any $t>0$,
\begin{eqnarray*}\label{def-e107a}
C(F,t):=\limsup\overline{F}(x-t)(\overline{F}(x))^{-1}<\infty.
\end{eqnarray*}
The inclusion ${\cal{OS}}\subset {\cal{OL}}$ is well-known.\\

The following proposition gives a positive answer to the previous question and has important implications concerning Theorem 19 of Beck et al. \cite{BBS13}.
At the same time it provides one of the four examples required to prove Theorem \ref{thm102}.

\begin{pron}\label{pron102}
There exists a distribution $F$ such that $F\in\mathcal{DK}^c\setminus \cal{OL}$ and $F\notin\mathcal{DK}_1$, thus $F\notin\J$, but $F^{I}\in(\mathcal{J}\cap\mathcal{DK}^c)\setminus({\cal{L}\cup\cal{D}})$, $F^I\notin\mathcal{DK}_1$, and $F^I$ is not weakly tail equivalent to a distribution in $\cal{S}$.
\end{pron}

We prove Proposition \ref{pron102}  in Section 4.

\section{Proof of Theorem \ref{thm101}}
\setcounter{equation}{0}\setcounter{Corol}{0}\setcounter{lemma}{0}\setcounter{re}{0}\setcounter{exam}{0}


We prove that the second inclusion in the theorem holds. Suppose that $F \in \mathcal{OS}\backslash \mathcal{DK}^c$.
Then there exists some $\lambda>0$ and a sequence of positive numbers $\{x_n,n\ge1\}$ such that $x_n>2x_{n-1}$  and
\begin{eqnarray}
\overline{F}(x_n)\le \exp\{-\lambda x_n\}
\label{301}
\end{eqnarray}
for every $n \in \mathbb{N}$. Since $F \in \mathcal{OS}$, there exist two constants
$2\le C^{*}(F)=\limsup\frac{\overline {F^{*2}}(x)}{\overline F(x)}<\infty$ and $y_0$ large enough such that
\begin{eqnarray*}
\overline {F^{*2}(y)}\le 2C^{*}(F) \overline {F}(y)
\end{eqnarray*}
for all $y\ge y_0$. Take any $y\ge y_0$, then
\begin{eqnarray*}
\overline{F}(y/2)&\le& \sqrt{ \mathbb{P}(S_2>y)}\le \sqrt{2C^{*}(F)\overline{F}(y)}.
\end{eqnarray*}
Iterating, for any positive integer $m$, we get
\begin{eqnarray}
\overline{F}(2^{-m}y)\le (2C^{*}(F))^{2^{-1}+\cdot\cdot\cdot+2^{-m}}(\overline{F}(y))^{2^{-m}}
\le 2C^{*}(F)(\overline{F}(y))^{2^{-m}},
\label{302}
\end{eqnarray}
as long as $2^{-m+1}y \ge y_0$. Without loss of generality,
we assume that $x_1\ge y_0$.
For any $x>x_1$, there exists a positive integer $n=:n(x)\ge2$ such that $x_{n-1}<x\le x_{n}$. Further, there exists a positive integer $m=:m(x)$ such that
$$
\max\{x_{n-1},2^{-m}x_n\}< x \le 2^{-m+1}x_n.
$$
If $2^{-m}x_n\ge x_{n-1}\ge y_0$, then by (\ref{302}) and (\ref{301}), we have
$$
\overline{F}(x) \le \overline{F}(2^{-m}x_n)\le 2C^{*}(F)(\overline{F}(x_n))^{2^{-m}}\le2C^{*}(F)\exp\{-\lambda 2^{-m} x_n\}\le2C^{*}(F)\exp\{-2^{-1}\lambda x\};
$$
if $2^{-m}x_n<x_{n-1}$. Thus, by (\ref{301}), we obtain that
$$
\overline{F}(x) \le \overline{F}(x_{n-1})\le\exp\{-\lambda x_{n-1}\}\le\exp\{-2^{-1}\lambda 2^{-m+1}x_{n}\}
\le\exp\{-2^{-1}\lambda x\},
$$
so $F$ is light-tailed and the claim therefore follows.\hfill$\Box$\\

The following example, which was introduced by C. M. Goldie (see Example 4.1 of \cite{K88}), shows that the inclusion
$\mathcal{OS}\cap\mathcal{K}\subset\mathcal{DK}^c$ is proper.

\begin{exam}\label{exam300}
Take $x_0=0,\ x_n=\sum_{k=0}^n k^{k-2},\ n\ge2$. Define a distribution $F$ such that
\begin{eqnarray*}
\overline{F}(x)=
\sum_{n=1}^\infty n^{-n}\textbf{\emph{1}}(x\in[x_{n-1},x_n))
\end{eqnarray*}
for $x \ge 0$. Then $F\in\mathcal{DK}^c$, but $F\notin\mathcal{OS}\cap\mathcal{K}$.
\end{exam}
\noindent{\bf Proof} \ For $n$ large enough, we have
$$\frac{\overline{F}(x_n-1)}{\overline{F}(x_n)}=\frac{\overline{F}(2^{-1}x_n)}{\overline{F}(x_n)}\ge n+1\rightarrow\infty,$$
as $n\rightarrow\infty$. Thus $F\notin\cal{OL}\supset\mathcal{OS}\cap\mathcal{K}$. But when $x\in[x_{n-1},x_n)$ and $\delta>1$, from
\begin{eqnarray*}
x^{\delta}\overline{F}(x)\ge x_{n-1}^{\delta}\overline{F}(x_{n})\ge(n-1)^{\delta(n-3)}n^{-n}\rightarrow\infty
\end{eqnarray*}
as $n\rightarrow\infty$, we obtain $F\in\mathcal{DK}_1\subset\mathcal{DK}^c$.\hfill$\Box$\\

\section{Proof of Theorem \ref{thm102}}
\setcounter{equation}{0}\setcounter{Corol}{0}\setcounter{lemma}{0}\setcounter{re}{0}\setcounter{exam}{0}

The first example shows that $(\J\cap \mathcal{K})\backslash (\mathcal{L}\cup \mathcal{D})$ contains distributions which are not in $\mathcal{DK}_1$
and which are weakly tail equivalent to a distribution in $\mathcal{S}$.

\begin{exam}\label{exam201}
Assume that $F_1\in\mathcal{S}$ is continuous with all (polynomial) moments finite and let $y_0\ge0$ and $a>1$ be two constants such that $a\overline{F_1}(y_0)\le1$. For example,
we can take $\overline{F_1}(x)=e^{-\sqrt{x}},\ x\in [0,\infty)$ and $y_0=(\ln a)^2$. Define the distribution $F$ by
\begin{eqnarray}\label{exam-e101}
\overline{F}(x)&=&\overline{F_1}(x)\textbf{\emph{1}}(x<x_1)+
\sum_{i=1}^{\infty}\Big(\overline{F_1}(x_i)\textbf{\emph{1}}(x_i\le x<y_i)+\overline{F_1}(x)\textbf{\emph{1}}(y_i\le x<x_{i+1})\Big),
\end{eqnarray}
where $\{x_i,i\ge1\}$ and $\{y_i, i\ge 1\}$ are two sequences of positive constants satisfying $x_i<y_i<x_{i+1}$ and $\overline{F_1}(x_i)=a\overline{F}_1(y_i),\ i\ge1$. Then $F\in(\J\cap\mathcal{DK}^c)\setminus(\cal{L}\cup\cal{D})$ and $F\notin\mathcal{DK}_1$. Obviously, $\overline{F}(x)\approx\overline{F_1}(x)$.
\end{exam}

\noindent{\bf Proof} \ It is easy to see that  $\overline{F_1}(x)\le\overline F(x)\le a\overline{F_1}(x)$. Since $F_1\in\mathcal{S}=\mathcal{L}\cap
\J$ and using the fact that $\J$ is closed under weak tail equivalence (see
Proposition 8 of \cite{BBS13}), we obtain that $F\in\mathcal{J}$. Next it is easy to verify that $F\in\mathcal{K}\setminus\mathcal{DK}_1$
and hence $F\notin\mathcal{D}$, and since $\lim_{n \to \infty} (\overline F(y_n))^{-1}\overline F(y_n-1)=a$ we get $F\notin\mathcal{L}.$
\hfill$\Box$\\

Next, we provide an example in $(\J \cap \mathcal{K})\backslash (\mathcal{L}\cup \mathcal{D})$ which is in $\mathcal{DK}_1$ and which is
weakly tail equivalent to a distribution in $\mathcal{S}$. To this end, we first construct a
distribution $F_1$ belonging to the class $(\mathcal{S}\cap \mathcal{DK}_1)\setminus\mathcal{D}$ and then show that $F$ defined as in the previous example
has the required properties.

\begin{exam}\label{exam202}
Choose any constants $\alpha\in
(0,1)$, $\beta\in(\alpha,2\alpha)$ and $x_1>2^{\alpha(\beta-\alpha)^{-1}}$. For all
integers $n\geq1$, let $x_{n+1}=x_n^{\beta{\alpha}^{-1}}$. Clearly,
$x_{n+1}>2x_n$ and $x_n\to\infty$ as $n\to\infty$. Now, define the
distribution $F_1$ as follows£º
\begin{eqnarray}\label{exam-e201}
\overline{F_1}(x)&=&(x_1^{-1}(x_1^{-\alpha}-1)x+1)\textbf{\emph{1}}(0\leq x< x_1)\nonumber\\
\ \ \ \ \ \ \ \ \ \ \ \ \ \ \ \ \ &+&\sum\limits_{n=1}^{\infty}(x_n^{-\alpha}+(x_n^{-\beta}-x_n^{-\alpha})(x_{n+1}-x_n)^{-1}(x-x_n))
\textbf{\emph{1}}(x_n\leq x<x_{n+1}),
\end{eqnarray}
where $x\in [0,\infty)$. Then $F_1$ has an infinite mean and belongs to the class $(\mathcal{S}\cap \mathcal{DK}_1)\setminus\mathcal{D}$.

Further, in the above example, we take $a=2$ and $y_n=2^{-1}(x_{n+1}+x_n)+2^{-1}(x_n^{\beta-\alpha}-1)^{-1}(x_{n+1}-x_n)$ such that $x_n<y_n<x_{n+1}$ and $2\overline{F_1}(y_n)=\overline{F_1}(x_n)$ for all $n\ge1$. Let $F$ be a distribution such that
\begin{eqnarray}\label{exam-e202}
\overline{F}(x)&=&\overline{F_1}(x)\textbf{\emph{1}}(x<x_1)+
\sum_{n=1}^{\infty}\bigg(\overline{F_1}(x_n)\textbf{\emph{1}}(x_n\le x<y_n)\nonumber\\
& &+\overline{F_1}(x)\textbf{\emph{1}}(y_n\le x<x_{n+1})\bigg),\ x\in(-\infty,\infty),
\end{eqnarray}
then $F\in(\J\cap\mathcal{DK}_1)\setminus(\cal{L}\cup\cal{D})$ and $\overline{F}(x)\approx\overline{F_1}(x)$.
\end{exam}

\noindent{\bf Proof} \ The distribution $F_1$ clearly has an infinite mean. It belongs to the class $\mathcal{DK}_1\setminus\mathcal{D}$ due to the following two facts
\begin{eqnarray*}
(\overline{F_1}(x_{n+1}))^{-1}\overline{F_1}(2^{-1}x_{n+1})\ge2^{-1}(1+x_n^{\beta-\alpha})\to \infty,~n\rightarrow\infty,
\end{eqnarray*}
where used the elementary inequality $\frac{a-c}{b-c}\le\frac{a}{b}$ for $b\ge a>0$ and $c\ge0$, and for all $x\in[x_n,x_{n+1})$, $n\ge1$
\begin{eqnarray*}
x\overline{F_1}(x)\ge x_n\overline{F_1}(x_n)\sim x_n^{1-\alpha}\to \infty,~n\rightarrow\infty.
\end{eqnarray*}

Next we prove $F_1\in\mathcal{S}$. By
\begin{eqnarray*}
\overline {F_1^{*2}}(x)&=&2\overline {F_1}(x)-\overline {F_1}^2(2^{-1}x)+2\int_{2^{-1}x}^{x}\overline {F_1}(x-y)F_1(dy),
\end{eqnarray*}
and
\begin{eqnarray*}
\liminf_{x\to\infty}(\overline {F_1}(x))^{-1}\overline {F_1^{*2}}(x)=2,
\end{eqnarray*}
we only need to prove
\begin{eqnarray*}
H(x):=(\overline {F_1}(x))^{-1}\int_{2^{-1}x}^{x}\overline {F_1}(x-y)F_1(dy)\to0,~x\rightarrow\infty.
\end{eqnarray*}
To this end, we estimate $H(x)$ in the  two cases $x_n\le x <2x_n$ and $2x_n\le x <x_{n+1},\ n\ge1$.
When $x\in[x_n,2x_n)$, then $2^{-1}x\le y \le x$ implies $2^{-1}x_{n}\le y \le x_n$, $n\ge1$. Thus by (\ref{exam-e201}) we have
\begin{eqnarray}\label{1010}
H(x)&\leq&(\overline {F_1}(2x_n))^{-1}(x_n-x_{n-1})^{-1}x_{n-1}^{-\alpha}\lo(\int_{0}^{x_{n-1}}+\int_{x_{n-1}}^{x_n}\ro)\overline {F_1}(y)dy\nonumber\\
&\le&(\overline {F_1}(2x_n))^{-1}(x_n-x_{n-1})^{-1}x_{n-1}^{-\alpha} \lo(x_{n-1}+\overline {F_1}(x_{n-1})(x_n-x_{n-1})\ro)\nonumber\\
&\sim&x_n^{-(1-\alpha)(1-\beta^{-1}\alpha)}+x_n^{-\alpha(2\beta^{-1}\alpha-1)}\to0,~n\rightarrow\infty.
\end{eqnarray}
When $x\in[2x_n,x_{n+1})$, then $2^{-1}x\le y \le x$ implies $x_n\le y \le x_{n+1}$, $n\ge1$. Thus by (\ref{exam-e201}) we obtain
\begin{eqnarray}\label{1011}
H(x)&\leq&(\overline {F_1}(x_{n+1}))^{-1}(x_{n+1}-x_n)^{-1}x_n^{-\alpha}\lo(\int_{0}^{x_n}+\int_{x_n}^{2^{-1}x_{n+1}}\ro)\overline {F_1}(y)dy\nonumber\\
&\le&(\overline {F_1}(x_{n+1}))^{-1}(x_{n+1}-x_n)^{-1}x_n^{-\alpha}\lo(x_n+\overline {F_1}(x_{n})(2^{-1}x_{n+1}-x_n)\ro)\nonumber\\
&\sim&x_n^{-(\beta-\alpha)(\alpha^{-1}-1)}+2^{-1}x_n^{-(2\alpha-\beta)}\to0,~n\rightarrow\infty.
\end{eqnarray}
By (\ref{1010}) and (\ref{1011}) we get $F_1\in\mathcal{S}$.

Now, we prove that $F\in(\J\cap\mathcal{DK}_1)\setminus(\cal{L}\cup\cal{D})$. By (\ref{exam-e202}) we easily get
\begin{eqnarray}\label{2010}
\overline{F_1}(x)\le \overline{F}(x) \le 2\overline{F_1}(x),
\end{eqnarray}
that is $\overline{F}(x)\approx\overline{F_1}(x)$. Then by (\ref{2010}) and $F_1\in\mathcal{S}\subset\J$ we have $F\in\mathcal{J}$. Next by $F_1\in\mathcal{DK}_1\setminus\mathcal{D}$ we immediately get $F\in\mathcal{DK}_1\setminus\mathcal{D}$.
Finally, $F\notin\mathcal{L}$ follows from
\begin{eqnarray*}
(\overline{F}(y_n))^{-1}\overline{F}(y_n-1)=2,
\end{eqnarray*}
for all $n\ge1$.\hfill$\Box$\\



The following example shows that there exist distributions in  $(\J\cap \mathcal{DK}_1)\setminus(\cal{L}\cup\cal{D})$ which are not weakly tail equivalent
to a distributon in $\mathcal{S}$.

\begin{exam}\label{exam401}
Let $m\geq 1$ be an any integer. Choose any constants $\alpha\in
(2+3m^{-1},\infty)$ and $x_1>4^{\alpha}$. For all
integers $n\geq1$, let $x_{n+1}=x_n^{1+{\alpha}^{-1}}$. Clearly,
$x_{n+1}>4x_n$ and $x_n\to\infty$ as $n\to\infty$. Now, define the
distribution $F$ as follows£º
\begin{eqnarray}\label{exam-e401}
&&\overline{F}(x)=(x_1^{-1}(x_1^{-\alpha}-1)x+1)\textbf{\emph{1}}(0\leq x< x_1)\nonumber\\
&+&\sum\limits_{n=1}^{\infty}((x_n^{-\alpha}+(x_n^{-\alpha-2}-x_n^{-\alpha-1})(x-x_n))\textbf{\emph{1}}(x_n\leq x<2x_n)
+x_n^{-\alpha-1}\textbf{\emph{1}}(2x_n\leq x< x_{n+1})),
\end{eqnarray}
where $x\in [0,\infty)$. Further, let $m \in \N$ and $$\overline{G_m}(x)=(\overline{F}(x))^m:=\overline{F}^m(x),~x\in(-\infty,\infty).$$
Then $G_m\in(\J\cap \mathcal{DK}_1)\setminus(\cal{L}\cup\cal{D})$ with finite mean and $G_m$ is not weakly tail equivalent to a distribution in $\cal{S}$.
\end{exam}

\begin{re} \label{remark401} This example shows that there are many such distributions since $m$, $\alpha$ and $x_1$ are arbitrary. Further, based on each of the above
distributions and using the method of Example 2.1, we can construct new distributions in the class $\mathcal{DK}_1\setminus(\cal{L}\cup\cal{D})$ which are not weakly
tail equivalent to a distribution in $\mathcal{S}$.
\end{re}

\noindent{\bf Proof} \ According to Proposition 12 b) in \cite{BBS13} and the following Lemma \ref{lemma401}, we need to prove
the above conclusion only for $m=1$, i.e.~$G_1=F$. By (\ref{exam-e401}), it is easy to see that when $x\geq x_1,$ then
\begin{eqnarray}
x^{-\alpha-1}\leq\overline{F}(x)\leq 2^{\alpha}x^{-\alpha}.\label{e101-2}
\end{eqnarray}
Thus $F\in\mathcal{DK}_1$. Moreover, using (\ref{exam-e401}) and (\ref{e101-2}), we see that the distribution $F$ has a finite mean which we denote by $\mu$. In fact, we get
\begin{eqnarray}
\int_0^{\infty}y^{4}\overline{F}(y)dy<\infty.\label{e101-0}
\end{eqnarray}
\vspace{0.2cm}

Observe that $F\notin{\cal{L}\cup\cal{D}}$ by (\ref{exam-e401}) and the following facts
\begin{eqnarray*}
\lo(\overline{F}(2x_n)\ro)^{-1}\overline{F}(2x_n-1)&=&2-x_n^{-1}
\to2,~~n\to\infty.
\end{eqnarray*}
and
\begin{eqnarray*}
\overline{F}(2x_n)(\overline{F}(x_n))^{-1}=x_n^{-1}
\to0,~~n\to\infty.
\end{eqnarray*}
\vspace{0.2cm}

Next we prove $F\in \J.$ For $x > 2K$, we obtain

\begin{eqnarray}
&&B(x):=P(X_{2,2}\leq K|X_1+X_2\ge x)\nonumber\\
&=&\lo(\frac{1}{2}\overline {F}^2(\frac{x}{2})+\int_{0}^{\frac{x}{2}}\overline {F}(x-y)dF(y)\ro)^{-1}\int_{0}^{K}\overline {F}(x-y)dF(y)\nonumber\\
&=&1-\lo(\frac{1}{2}\overline {F}^2(\frac{x}{2})+\int_{0}^{\frac{x}{2}}\overline {F}(x-y)dF(y)\ro)^{-1}\lo(\int_{K}^{\frac{x}{2}}\overline {F}(x-y)dF(y)+\frac{1}{2}\overline {F}^2(\frac{x}{2})\ro)\nonumber\\
&\geq&1-\lo(\int_{0}^{\frac{x}{2}}\overline {F}(x-y)dF(y)\ro)^{-1}\lo
(\int_{K}^{\frac{x}{2}}\overline {F}(x-y)dF(y)\ro)
-\overline {F}^2(\frac{x}{2})\Big(2\overline {F}(x)F(\frac{x}{2})\Big)^{-1}\nonumber\\
&:=&1-B_1(x)-B_2(x).\label{e101-3}
\end{eqnarray}

By (\ref{e101-2}),  we have for $x \ge 2x_1$
\begin{eqnarray}
B_2(x)&\leq&(F(\frac{x}{2}))^{-1}2^{4\alpha-1}x^{1-\alpha}\to0.
\label{e101-4}
\end{eqnarray}

Next we estimate $B_1(x)$ in each of the five cases $x_n\leq x< x_n+K$, $x_n+K\leq x< 2x_n$,
$2x_n\leq x< 2x_n+K$, $2x_n+K\leq x< 4x_n$ and $4x_n\leq x< x_{n+1}$.

When $x\in[x_n, x_n+K)$, then $K\le y\le 2^{-1}x$ implies $2x_{n-1}\le x-y\le x_n, n\geq2$. Thus, by (\ref{exam-e401}), we obtain
\begin{eqnarray}
B_1(x)&\leq&\lo(\int_{0}^{K}\overline {F}(x-y)dF(y)\ro)^{-1}\int_{K}^{\frac{x}{2}}\overline {F}(x-y)dF(y)\nonumber\\
&\leq&\lo(\overline {F}(x_n+K)F(K)\ro)^{-1}\overline {F}(K)\overline {F}(2x_{n-1})\nonumber\\
&\sim&\lo(F(K)\ro)^{-1}\overline {F}(K)\to 0,~K\rightarrow\infty.
\label{e101-5}
\end{eqnarray}

When $x\in[x_n+K, 2x_n), n\geq2,$ by (\ref{exam-e401}),  we have,
\begin{eqnarray*}
B_1(x)&\leq&\lo(\int_{0}^{K}\overline {F}(x-y)dF(y)\ro)^{-1}{\lo(\int_{K}^{x-x_n}+\int_{x-x_n}^{\frac{x}{2}}\ro)}\overline {F}(x-y)dF(y)\nonumber\\
&:=&B_{11}(x)+B_{12}(x).
\end{eqnarray*}
Note that $x_n\le x-y\le2x_n$ for $K\le y\le x-x_n, n\ge 2$, so by (\ref{exam-e401}), (\ref{e101-2}) and (\ref{e101-0}), we have
\begin{eqnarray}
B_{11}(x)&\leq&\lo(\overline {F}(x)F(K)\ro)^{-1}\int_{K}^{x-x_n}(\overline F(x)+x_n^{-\alpha-1}y)dF(y)\nonumber\\
&\leq&\lo(F(K)\ro)^{-1}\int_{K}^{\infty}(1+y)dF(y)\to 0,~K\rightarrow\infty.
\label{e101-6}
\end{eqnarray}
Now we deal with $B_{12}(x)$ in the two cases $x_n+K\leq
x< \frac{3}{2}x_n$ and $\frac{3}{2}x_n\leq x< 2x_n$.
When $x\in[x_n+K,\frac{3}{2}x_n)$, then $2x_{n-1}\le x-y\le x_n$ for $x-x_n\le y\le \frac{x}{2}, n\geq2$, so by (\ref{exam-e401}), we have
\begin{eqnarray}
B_{12}(x)&\leq&\lo(\overline {F}(\frac{3}{2}x_n)F(K)\ro)^{-1}\lo(\overline F(2x_{n-1})\overline F(K)\ro)\nonumber\\
&\sim&\lo(F(K)\ro)^{-1}2\overline F(K)\to 0,~K\rightarrow\infty.
\label{e101-7}
\end{eqnarray}
When $x\in[\frac{3}{2}x_n,2x_n)$, then $2x_{n-1}\le y\le x_n, n\geq2$ and,  by (\ref{exam-e401}),  we get
\begin{eqnarray}
B_{12}(x)=0.
\label{e101-8}
\end{eqnarray}

When $x\in[2x_n,2x_n+K)$, if $K\le y\le x_n$, then $x_n\le 2x_n-y\le 2x_n$; if $x_n\le y\le x_n+2^{-1}K$, then $x_n-2^{-1}K\le 2x_n-y\le x_n, n\geq2$. Thus, by (\ref{exam-e401}), (\ref{e101-2}) and (\ref{e101-0}), we have
\begin{eqnarray}
B_1(x)&\leq&\lo(\overline {F}(x)F(K)\ro)^{-1}\lo(\int_{K}^{x_n}+\int_{x_n}^{x_n+\frac{K}{2}}\ro)\overline {F}(2x_n-y)dF(y)\nonumber\\
&\leq&\lo(F(K)\ro)^{-1}\lo(\int_{K}^{x_n}(1+y)dF(y){+\int_{x_n}^{x_n+\frac{K}{2}} x_ndF(y)}\ro)\nonumber\\
&\leq&\lo(F(K)\ro)^{-1}\lo(\int_{K}^{\infty}(1+y)dF(y){+\int_{K}^{\infty} ydF(y)}\ro)\to 0,~~~K\rightarrow\infty.
\label{e101-9}
\end{eqnarray}

When $x\in[2x_n+K,4x_n)$, if $K\le y\le x-2x_n$, then $2x_n\le x-y\le 4x_n$; if $x-2x_n\le y\le 2^{-1}x$, then $x_n\le x-y\le 2x_n, n\geq2$. Thus, by (\ref{exam-e401}), (\ref{e101-2}) and (\ref{e101-0}), we have
\begin{eqnarray}
B_1(x)&\leq&\lo(\int_{0}^{K}\overline {F}(x-y)dF(y)\ro)^{-1}\lo(\int_{K}^{x-2x_n}+\int_{x-2x_n}^{\frac{x}{2}}\ro)\overline {F}(x-y)dF(y)\nonumber\\
&\leq&\lo(\overline {F}(2x_n)F(K)\ro)^{-1}\Big(\overline F(2x_n)\overline F(K)\nonumber\\
&&+\int_{x-2x_n}^{\frac{x}{2}}\lo(x_n^{-\alpha}+(x_n^{-\alpha-2}-x_n^{-\alpha-1})(x-x_n-y)\ro)dF(y)\Big)\nonumber\\
&\leq&\lo(F(K)\ro)^{-1}\lo(\overline F(K)+\int_{K}^{\infty}(1+y)dF(y)\ro)\to 0,~~~K\rightarrow\infty.
\label{e101-10}
\end{eqnarray}

When $x\in[4x_n,x_{n+1})$, if $0\le y\le 2^{-1}x$, then $2x_n\le x-y\le x_{n+1}, n\geq2$, so by (\ref{exam-e401}), we have
\begin{eqnarray}
B_1(x)&\leq&{\lo(F(K)\ro)^{-1}\overline F(K)}\to 0,~K\rightarrow\infty.
\label{e101-11}
\end{eqnarray}

By (\ref{e101-3})-(\ref{e101-11}), we get $F\in \J$.
\vspace{0.2cm}

Next, we prove that $F$ is not weakly tail equivalent to a distribution in $\mathcal{S}$. To see this, we state the following lemma.

\begin{lemma}\label{lemma401}
Assume that the distribution $F$ satisfies
\begin{eqnarray}
\limsup_{t\rightarrow\infty} C(F,t)=\limsup_{t\rightarrow\infty}\limsup\overline{F}(x-t)(\overline{F}(x))^{-1}=\infty.\label{2670}
\end{eqnarray}
Then $F$ is not weakly equivalent to any long-tailed distribution.
\end{lemma}

\noindent{\bf Proof} We assume there exists a distribution $F_1\in\mathcal{L}$ and $\overline{F_1}(x)\approx \overline{F}(x)$. Then there are two constants $0<C_1\leq C_2<\infty$ such that
\begin{eqnarray}\label{1000}
C_1= \liminf(\overline{F}(x))^{-1}\overline{F_1}(x)\leq \limsup(\overline{F}(x))^{-1}\overline{F_1}(x)= C_2.
\end{eqnarray}
By $F_1\in\mathcal{L}$ and (\ref{1000}) for any $0<t<\infty$ we have
\begin{eqnarray}\label{2000}
(\overline{F}(x))^{-1}\overline{F}(x-t)
\leq(C_1\overline{F_1}(x))^{-1}C_2\overline{F_1}(x-t)\sim&C_1^{-1}C_2.
\end{eqnarray}
Obviously, (\ref{2000}) contradicts  (\ref{2670}). Hence the conclusion of the lemma holds.\hfill$\Box$\\

In Example \ref{exam401}, we have
\begin{eqnarray*}
(\overline{F}(2x_n))^{-1}\overline{F}(2x_n-t)
=1+t-tx_n^{-1}\rightarrow 1+t,~~~n\to\infty,
\end{eqnarray*}
that is (\ref{2670}) holds. Thus, by Lemma \ref{lemma401}, $F$ is not weakly tail equivalent to a distribution in $\mathcal{S}$.\hfill$\Box$\\

The remaining statement of Theorem \ref{thm102} will be proved in the next section.

\section{Proof of Proposition \ref{pron102}}
\setcounter{equation}{0}\setcounter{Corol}{0}\setcounter{lemma}{0}\setcounter{re}{0}\setcounter{exam}{0}

Next, we provide an example proving Proposition \ref{pron102} and which at the same time shows the remaining claim of Theorem \ref{thm102} that
there exists a distribution in $(\J\cap\mathcal{DK}^c)\setminus({\cal{L}\cup\cal{D}})$ which is not in  $\mathcal{DK}_1$ and which is not
weakly tail equivalent to a distribution in $\mathcal{S}$.

\begin{exam}\label{exam501}
Define the distribution $F$ as follows
\begin{eqnarray}\label{exam-e102}
\overline{F}(x)=8^{-1}\textbf{\emph{1}}(0\le x<4)
+\sum\limits_{n=2}^{\infty}(2^{\frac{-n-n^2}{2}}-2^{\frac{-3n-n^2}{2}})\textbf{\emph{1}}(2^n\le x<2^{n+1})),
\end{eqnarray}
where $x\in [0,\infty)$. Then $F\in\mathcal{DK}^c\setminus \cal{OL}$, thus $F\notin\J$, but $F^{I}\in(\J\cap\mathcal{DK}^c)\setminus({\cal{L}\cup\cal{D}})$, $F^{I}\notin\mathcal{DK}_1$, and $F^{I}$ is not weakly tail equivalent to a distribution in $\cal{S}$.
\end{exam}

\noindent{\bf Proof} \ First, it is easy to check that
\begin{eqnarray*}
\overline {F}(2^{n}-1)(\overline {F}(2^{n}))^{-1}
&\sim&2^n\to\infty,
~~~n\rightarrow\infty,
\end{eqnarray*}
so $F\not\in{\cal{OL}}$. Since $\J\subset\cal{OL}$, we get $F\not\in \J$. It is easy to verify that $F \in \mathcal{DK}^c\setminus\mathcal{DK}_1$. In addition, since
\begin{eqnarray*}
\int_{0}^{\infty}\overline {F}(y)dy=2^{-1}+\sum_{n=2}^{\infty}(2^{\frac{n-n^2}{2}}-2^{\frac{-n-n^2}{2}})=1,
\end{eqnarray*}
the distribution $F$ has a finite mean $\mu=1$. Further, we obtain
\begin{eqnarray}\label{newlabel}
&&\overline {F^I}(x)=(1-\frac{x}{8})\mathrm{\bf{1}}(0\le x<4)\nonumber\\
&+&\sum_{n=2}^{\infty}\lo(2^{\frac{n-n^2}{2}}+(2^{\frac{-3n-n^2}{2}}-2^{\frac{-n-n^2}{2}})(x-2^n)\ro)\mathrm{\bf{1}}(2^n\le x<2^{n+1}), ~x\in [0,\infty).\label{e102-1}
\end{eqnarray}

It is easy to see that $F^I \in \mathcal{DK}^c\setminus\mathcal{DK}_1$ and $F^I\not\in{\cal{L}\cup\cal{D}}$ by (\ref{e102-1}) and the following facts
\begin{eqnarray*}
\overline {F^I}(2^{n+1}-1)(\overline {F^I}(2^{n+1}))^{-1}&=&2+2^{-n}\to2, ~n\rightarrow
\infty
\end{eqnarray*}
and
\begin{eqnarray*}
\overline{F^I}(2^{n+1})(\overline{F^I}(2^n))^{-1}=2^{-n}\to 0, ~n\rightarrow
\infty.
\end{eqnarray*}

Next we prove $F^I\in \J.$ Let $\{X_i,i\geq 1\}$ be a sequence of i.i.d. random variables with common distribution $F^I$. As before, $X_{n,k}$
denotes the $k$-th largest random variable in the sequence $\{X_i,1\leq i\leq n\},1\leq k\leq n$. For $x>2K>0$ we obtain
\begin{eqnarray}
&&B(x):=\P(X_{2,2}\leq K|S_2\ge x)\nonumber\\
&=&\lo(\frac{1}{2}\overline {F^I}^2(\frac{x}{2})+\int_{0}^{\frac{x}{2}}\overline {F^I}(x-y)dF^I(y)\ro)^{-1}\int_{0}^{K}\overline {F^I}(x-y)dF^I(y)\nonumber\\
&=&1-\lo(\frac{1}{2}\overline {F^I}^2(\frac{x}{2})+\int_{0}^{\frac{x}{2}}\overline {F^I}(x-y)dF^I(y)\ro)^{-1}\lo(\int_{K}^{\frac{x}{2}}\overline {F^I}(x-y)dF^I(y)+\frac{1}{2}\overline {F^I}^2(\frac{x}{2})\ro)\nonumber\\
&\geq&1-\lo(\int_{0}^{\frac{x}{2}}\overline {F^I}(x-y)dF^I(y)\ro)^{-1}\lo(\int_{K}^{\frac{x}{2}}\overline {F^I}(x-y)dF^I(y)\ro)-\frac{\overline {F^I}^2(\frac{x}{2})}{2\overline {F^I}(x)F^I(\frac{x}{2})}\nonumber\\
&=:&1-B_1(x)-B_2(x).\label{e102-2}
\end{eqnarray}

For all $x\in[2^n,2^{n+1}),n\geq2,$ by (\ref{e102-1}), we have
\begin{eqnarray}
B_2(x)\leq(2\overline {F^I}(2^{n+1})F^I(2^{n-1}))^{-1}\overline {F^I}^2(2^{n-1})
=(F^I(2^{n-1}))^{-1}2^{\frac{7n-6-n^2}{2}}\to0,
~n\rightarrow\infty.
\label{e102-3}
\end{eqnarray}

Now, we estimate $B_1(x)$ in the two cases $2^n\leq
x< 2^n+K$ and $2^n+K\leq x< 2^{n+1}$.

When $x\in[2^n,2^{n}+K),n\geq2,$ by {(\ref{e102-1})} and $\int_{0}^{\infty}ydF^I(y)<\infty$, we have

\begin{eqnarray}
&&B_1(x)\leq{\lo(\int_{0}^{K}\overline {F^I}(x-y)dF^I(y)\ro)^{-1}\int_{K}^{\frac{x}{2}}\overline {F^I}(x-y)dF^I(y)}\nonumber\\
&\leq&\lo(\overline {F^I}(x)F^I(K)\ro)^{-1}\int_{K}^{\frac{x}{2}}2^{\frac{3n-2-n^2}{2}}+
(2^{\frac{2-n-n^2}{2}}-2^{\frac{n-n^2}{2}})(x-K-2^{n-1}-y+K)dF^I(y)\nonumber\\
&\leq&\lo(\overline {F^I}(x)F^I(K)\ro)^{-1}\int_{K}^{\frac{x}{2}}\lo(\overline {F^I}(x-K)+2^{\frac{n-n^2}{2}}y\ro)dF^I(y)\nonumber\\
&\leq&\lo(\overline {F^I}(2^n+K)F^I(K)\ro)^{-1}\lo(\overline {F^I}(2^n-K)\overline {F^I}(K)+\int_{K}^{\frac{x}{2}}2^{\frac{n-n^2}{2}}ydF^I(y)\ro)\nonumber\\
&=&O\lo(\overline {F^I}(K)+\int_{K}^{+\infty}ydF^I(y)\ro)\to 0,~~~K\rightarrow\infty.
\label{e102-4}
\end{eqnarray}
When $x\in[2^{n}+K,2^{n+1}),n\geq2,$ by {(\ref{e102-1})}, we have
\begin{eqnarray}
B_1(x)&\leq&{\lo(\int_{0}^{x-2^n}\overline {F^I}(x-y)dF^I(y)\ro)^{-1}\lo(\int_{K}^{x-2^n}+\int_{x-2^n}^{\frac{x}{2}}\ro)\overline {F^I}(x-y)dF^I(y)}\nonumber\\
&:=&B_{11}(x)+B_{12}(x).\label{e102-5}
\end{eqnarray}
By (\ref{e102-1}), (\ref{e102-5}) and $\int_{0}^{\infty}ydF^I(y)<\infty$,  we have
\begin{eqnarray}
B_{11}(x)&=&{\lo(\int_{0}^{x-2^n}\overline {F^I}(x-y)dF^I(y)\ro)^{-1}\lo(\int_{K}^{x-2^n}\overline {F^I}(x-y)dF^I(y)\ro)}\nonumber\\
&\leq&\lo(\overline {F^I}(x)F^I(x-2^n)\ro)^{-1}\lo(\int_{K}^{x-2^n}\lo(\overline {F^I}(x)+2^{\frac{-n-n^2}{2}}y\ro)dF^I(y)\ro)\nonumber\\
&\leq&\lo(F^I(K)\ro)^{-1}\lo(\int_{K}^{x-2^n}\lo(1+y\ro)dF^I(y)\ro)\to 0,~~~K\rightarrow\infty.
\label{e102-6}
\end{eqnarray}
Now we deal with $B_{12}(x)$ in the two cases $2^n+K\leq
x<3\times2^{n-1}$ and $3\times2^{n-1}\leq x< 2^{n+1}$.
When $x\in[2^n+K,3\times2^{n-1}),n\geq2,$ by (\ref{e102-1}), (\ref{e102-5}) and $\int_{0}^{\infty}ydF^I(y)<\infty$, we have
\begin{eqnarray}
B_{12}(x)&\leq&\lo(\overline {F^I}(x)F^I(x-2^n)\ro)^{-1}\lo(\int_{x-2^n}^{\frac{x}{2}}2^{\frac{3n-2-n^2}{2}}+
(2^{\frac{2-n-n^2}{2}}-2^{\frac{n-n^2}{2}})(x-2^{n-1}-y)dF^I(y)\ro)\nonumber\\
&\leq&\lo(\overline {F^I}(3\times2^{n-1})F^I(K)\ro)^{-1}\lo(\int_{K}^{\frac{x}{2}}2^{\frac{n-n^2}{2}}dF^I(y)+\int_{K}^{\frac{x}{2}}2^{\frac{n-n^2}{2}}ydF^I(y)\ro)\nonumber\\
&\leq&\lo(F^I(K)\ro)^{-1}\lo(2\overline {F^I}(K)+2\int_{K}^{+\infty}ydF^I(y)\ro)\to 0,~~~K\rightarrow\infty.
\label{e102-7}
\end{eqnarray}
When $x\in[3\times2^{n-1},2^{n+1}),n\geq2,$ by (\ref{e102-1}) and (\ref{e102-5}), we have
\begin{eqnarray}
B_{12}(x)&\leq&\lo(\overline {F^I}(x)F^I(x-2^n)\ro)^{-1}\overline {F^I}(2^{-1}x)\overline {F^I}(x-2^n)\nonumber\\
&\leq&\lo(\overline {F^I}(2^{n+1})F^I(2^{n-1})\ro)^{-1}(\overline {F^I}(2^{n-1}))^2\to0,~~~n\rightarrow\infty.
\label{e102-8}
\end{eqnarray}

Therefore, by (\ref{e102-2})-(\ref{e102-8}), we get
$F^I\in \J.$

Finally, by Lemma \ref{lemma401} and
\begin{eqnarray*}
(\overline{F^I}(2^{n+1}))^{-1}\overline{F^I}(2^{n+1}-t)
=1+t-2^{-n}t\sim1+t,~~~n\to\infty,
\end{eqnarray*}
we see that $F^I$ is not weakly tail equivalent to a distribution in $\cal{S}$.\hfill$\Box$\\

\begin{re}\label{remark501}
Observe that the example shows that $F^I \in \J$ does not imply that $F \in \J$. Conversely, $F \in \J$ does not imply $F^I \in \J$ either even if
we assume that $F$ has a finite first moment (otherwise $F^I$ is not defined). As an example,
we can take the example in Section 3.8 in \cite{FKZ13} for which $F \in \mathcal{S}$ and $F^I \notin \mathcal{S}$.  Since $\mathcal{S} \subset \J$ we
have $F \in \J$. Further, $\mathcal{S} \subset \mathcal{L}$ and $F \in \mathcal{L}$ implies $F^I \in \mathcal{L}$ by Lemma 2.26 in \cite{FKZ13}. Since
$\J \cap \mathcal{L}=\mathcal{S}$ we obtain $F^I \notin \J$.
\end{re}

\noindent\textbf{Acknowledgement.} The authors are very grateful to Sergey Foss (Heriot-Watt University, Edinburgh) for helpful discussions and suggestions.

\end{document}